\documentclass[final,1p,times]{elsarticle}

\usepackage{color}
\usepackage{caption}
\usepackage{amssymb}
\usepackage{geometry}
\usepackage{epstopdf}
\usepackage{subfigure}
\usepackage{longtable}
\usepackage{amsmath,amsfonts}
\usepackage{amsthm}
\usepackage{amssymb,version,graphicx,fancybox,mathrsfs,pifont,booktabs,mathtools}
\usepackage{url,hyperref,multirow}
\newcommand{\bs}[1]{\boldsymbol{#1}}
\newtheorem{thm}{\bf Theorem}[section]
\def \ri {{\rm i}}
\theoremstyle{plain}
\newtheorem{theorem}{Theorem}[section]
\newtheorem{lemma}[theorem]{Lemma}

\theoremstyle{definition}

\theoremstyle{remark}

\catcode`\@=11 \theoremstyle{plain}
\@addtoreset{equation}{section}   

\@addtoreset{figure}{section}
\renewcommand\thefigure{\thesection.\@arabic\c@figure}
\@addtoreset{table}{section}
\renewcommand\thetable{\thesection.\@arabic\c@table}

\begin{document}
\captionsetup[figure]{labelfont={bf},labelformat={default},labelsep=period,name={Fig.}}
\begin{frontmatter}

\title{High-order and Mass-conservative Regularized Implicit-explicit relaxation Runge-Kutta methods for the logarithmic Schr\"{o}dinger equation}

\footnote{{\small $^*$Corresponding author.

E-mail address: weiyb@buaa.edu.cn  (Y. Wei)}}

\author {Jingye Yan$^{1}$, Hong Zhang$^{2}$, Yabing Wei$^{3*}$, Xu Qian$^{2}$}
\address{$^1$ College of Mathematics and Physics, Wenzhou University, Wenzhou 325035, China\\
$^2$ Department of Mathematics, National University of Defense
Technology, Changsha 410073, China \\
$^3$ School of Mathematical Sciences, Jiangsu University, Zhenjiang 212013, China\\
}
\begin{abstract}
The non-differentiability of the singular nonlinearity (such as $f=\ln|u|^2$) at $u=0$ presents significant challenges in devising accurate and efficient numerical schemes for the logarithmic Schr\"{o}dinger equation (LogSE). To address this singularity, we propose an energy regularization technique for the LogSE. For the regularized model, we utilize Implicit-Explicit Relaxation Runge-Kutta methods, which are linearly implicit, high-order, and mass-conserving for temporal discretization, in conjunction with the Fourier pseudo-spectral method in space. Ultimately, numerical results are presented to validate the efficiency of the proposed methods.
\end{abstract}


\end{frontmatter}

\section{Introduction}
In this paper, we focus on the numerical solution for the logarithmic Schr\"{o}dinger equation (LogSE), which take the form:
\begin{equation}\label{logse}
\begin{cases}
\ri  \partial_t u(\bs x,t)+ \Delta u(\bs x,t)=\lambda uf(u), \quad   \bs x\in  \Omega, \;\; t>0,\\[2pt]
u(\bs x,0)=u_0(\bs x),\quad \bs x\in \bar \Omega,
\end{cases}
\end{equation}
where $\ri=\sqrt{-1},$ $f(u)=\ln(|u(\bs x,t)|^2),$ $\Omega\subset {\mathbb R}^d$, with $d\ge 1,$  represents a bounded domain with a smooth boundary. The function $u_0$ is a given initial condition with a regularity that will be specified later. The constant $\lambda\neq 0$ is a real number, where for $\lambda>0$, the solution exhibits repulsive or defocusing behavior, and for $\lambda<0$, the solutions demonstrate attractive or focusing interactions. The LogSE \eqref{logse} proposed as a model for nonlinear wave mechanics \cite{Bialynicki1976nonlinear}.

LogSE \eqref{logse} conserves both mass and energy as follows,
\begin{equation}\label{mass}
\begin{split}
M(t):&=\int_{\Omega}|u(\bs x, t)|^2 d \bs x \equiv \int_{\Omega}|u_0(\bs x)|^2 d \bs x=M(0),\\
E(t): & =\int_{\Omega}[|\nabla u(\bs x, t)|^2 d \bs x+\lambda F(|u(\bs x, t)|^2)] d \bs x \\
& \equiv \int_{\Omega}[|\nabla u_0(\bs x)|^2+\lambda F(|u_0(\bs x)|^2)] d \bs x=E(0),
\end{split}
\end{equation}
where
\begin{equation}\label{F}
F(\rho)=\int_0^\rho \ln (s) d s=\rho \ln \rho-\rho, \quad \rho =|u|^2.
\end{equation}

The nonlinearity term $f(u)=\ln(|u(\bs x,t)|^2)$ is non-differentiable at $u=0$. This characteristic introduces significant practical challenges and theoretical complexities, particularly when it comes to analyzing and solving the numerical solutions of the LogSE \eqref{logse}. Consequently, the existing literature on this topic is somewhat limited. For the Cauchy problem associated with the LogSE, Cazenave \cite{Cazenave1980} established a suitable functional framework.

Extensive literature exists on numerical solutions for the Schr\"odinger equation with smooth nonlinearity, employing methods such as the finite difference method \cite{Bao2019Regularized}, time-splitting method \cite{Bao2022Error,Bao2019Error}, and among others. However, for the LogSE \eqref{logse}, the available numerical schemes are more limited. This scarcity arises because these methods cannot be directly applied to the LogSE due to the non-differentiability of the nonlinearity at $u=0$. To circumvent this issue, Bao \cite{Bao2019Regularized,Bao2019Error} introduced a regularized logarithmic Schr\"odinger equation by substituting $f(u)$ with $f_{\varepsilon}(u^{\varepsilon})$, where $0<\varepsilon\ll 1$. They developed a semi-implicit finite difference method \cite{Bao2019Regularized}, which, although not conserving energy, provided a viable approach. Later, Bao \cite{Bao2022Error} introduced a different regularization strategy at the energy density level, substituting the energy density locally in the region $0<\rho<\varepsilon^2$ with a sequence of polynomials while keeping it unchanged in the region $\rho\geq\varepsilon^2$. In \cite{Bao2022Error,Bao2019Error}, Bao constructed regularized Lie-Trotter splitting and Strang splitting methods, which preserve mass conservation but are only first and second order, respectively, for solving the LogSE. Recently, Wang \cite{Wang2024Error} developed a nonregularized first-order implicit-explicit scheme for the LogSE, with the nonlinearity being explicit but only first order in time.

Over the past decades, a variety of relaxation Runge-Kutta (RRK) methods \cite{Ketcheson2019Relaxation} have been developed, offering explicit and mass conservative solutions. However, these explicit methods are plagued by stringent step size restrictions that can limit their practical applicability. Recently, Li \cite{Li2023Implicit} investigated implicit-explicit relaxation Runge-Kutta (IMEX RRK) methods for nonlinear stiff ordinary differential equations. The proposed methods are linearly implicit, can achieve arbitrarily high order accuracy, and are designed to preserve monotonicity.

The structure of the rest of this paper is as follows: Section 2 introduces the energy regularized logarithmic Schr\"{o}dinger equation (ERLogSE). In Section 3, IMEX RRK methods are applied to discretize the ERLogSE in time, combined with the Fourier pseudo-spectral method for spatial discretization. Section 4 presents several numerical experiments to demonstrate the efficiency and accuracy of the proposed numerical methods. Finally, Section 6 concludes the paper.
\section{The energy regularized LogSE }
In this context, the energy-regularization technique proposed by Bao et al. in \cite{Bao2022Error} is employed to address the singularity issue of the logarithmic term in the LogSE. This approach uses a polynomial approximation to smooth out the singularity at the origin for the energy density function $F(\rho)$, as depicted in equation \eqref{F}. The energy regularized method, as demonstrated in \cite{Bao2022Error}, yields superior performance compared to direct regularization of the logarithmic nonlinearity in the LogSE. The energy density $F(\rho)$ is approximated using a piecewise smooth function that incorporates a polynomial approximation near the origin as follows,
\begin{equation}\label{enerreg1}
F_{k}^{\varepsilon}(\rho)=F(\rho) \chi_{\{\rho \geq \varepsilon^{2}\}}+P_{k+1}^{\varepsilon}(\rho) \chi_{\{\rho<\varepsilon^{2}\}}, \quad k \geq 2,
\end{equation}
where $\chi_{_A}$ is the characteristic function of the set $A$ and
\begin{equation*}
P_{k+1}^{\varepsilon}(\rho)=\rho \ln \varepsilon^{2}-1-\sum \limits_{j=1}^{k} \frac{1}{j}(1-\frac{\rho}{\varepsilon^{2}})^{j}
\end{equation*}
 is a polynomial of degree $k+1$, which allows $F_{k}^{\varepsilon} \in C^{k}([0,+\infty)), F_{k}^{\varepsilon}(0)=F(0)=0$ (this satisfies the regularized energy to be well-defined on the whole space).
We can obtain
\begin{equation}\label{flose}
\begin{split}
f_{k}^{\varepsilon}(\rho)&=(F_{k}^{\varepsilon})^{\prime}(\rho)=\ln \rho \chi_{\{\rho \geq \varepsilon^{2}\}}+(P_{k+1}^{\varepsilon})^{\prime}(\rho) \chi_{\{\rho<\varepsilon^{2}\}}\\
&=\ln \rho \chi_{\{\rho \geq \varepsilon^{2}\}}+[\ln (\varepsilon^{2})-\frac{k+1}{k}(1-\frac{\rho}{\varepsilon^{2}})^{k}-\sum_{j=1}^{k-1} \frac{1}{j}(1-\frac{\rho}{\varepsilon^{2}})^{j}]\chi_{\{\rho<\varepsilon^{2}\}}, \quad \rho \geq 0.
\end{split}
\end{equation}
One can get the energy regularized logarithmic Schr\"odinger equation (ERLogSE) with a small parameter $0<\varepsilon\ll 1$ as follows
\begin{equation}\label{rlogse}
\begin{cases}
\ri  \partial_t u^{\varepsilon}(\bs x,t)+ \Delta u^{\varepsilon}(\bs x,t)=\lambda u^{\varepsilon} f^k_{\varepsilon}(u^{\varepsilon}), \quad   \bs x\in  \Omega, \;\; t>0,\\[2pt]
 u^{\varepsilon}(\bs x,0)=u_0(\bs x),\quad \bs x\in \bar \Omega,
\end{cases}
\end{equation}
which preserves the mass conservation law
\begin{equation}\label{rmass}
\begin{split}
M^{\varepsilon}(t):&=\int_{\Omega}|u^{\varepsilon}(\bs x, t)|^2 d \bs x \equiv \int_{\Omega}|u_0(\bs x)|^2 d \bs x=M(0),\\
\end{split}
\end{equation}
where
$f^k_\varepsilon$ defined in \eqref{flose} for ERLogSE.
\section{The construction of Regularized IMEX RRK methods}
\subsection{High order Regularized IMEX RRK methods}
 We rewrite ERLogSE \eqref{rlogse} as
 \begin{equation}\label{logse1}
 u^{\varepsilon}_t=g^{I}(u^\varepsilon)+g^{E}(u^\varepsilon),
 \end{equation}
 where $g^{I}(u^\varepsilon)=\ri \Delta u^\varepsilon,\, g^{E}(u^\varepsilon)=\ri \lambda {u}^{\varepsilon}f^k_\varepsilon(u^\varepsilon)$.

Applying the different orders IMEX RRK methods to the ERLogSE \eqref{logse1}.  We can get the
semi-discrete system as follows:
\begin{equation}\label{RRK}
\begin{cases}
u^\varepsilon_{i}=u^{\varepsilon,n}_{\gamma}+\tau\sum\limits^{i}_{j=1}a^{I}_{ij}g^{I}_{j}+\tau\sum\limits^{i-1}_{j=1}a^{E}_{ij}g^{E}_{j},\,\, i=1,\ldots,s,\\
u^{\varepsilon,n+1}_{\gamma}=u^{\varepsilon,n}_{\gamma}+\tau\gamma_{n}\sum\limits^{s}_{j=1}b^{I}_{j}g^{I}_{j}+\tau\gamma_{n}\sum\limits^{s}_{j=1}b^{E}_{j}g^{E}_{j},
\end{cases}
\end{equation}
where $g^{I}_{j}=g^{I}(t_n+c_j\tau,u^\varepsilon_{j}),\, g^{E}_{j}=g^{E}(t_n+c_j\tau,u^\varepsilon_{j}),\,\,j=1,\ldots,s$, and $u^{\varepsilon,n+1}_{\gamma}$ is the approximation at $\hat{t}_{n+1}=\hat{t}_{n}+\gamma_n\tau$. Define
\begin{equation*}
\begin{split}
&A^{I}= ({a^{I}_{ij}})_{s\times s},\,\, a^{I}_{ij}=0 \,\,\text{for} \,\,j>i,\\
&A^{E}= ({a^{E}_{ij}})_{s\times s},\,\, a^{I}_{ij}=0 \,\,\text{for} \,\,j\geq i,\\
&b^{I}=(b^{I}_{1},\ldots,b^{I}_{s})^{T},\,\,b^{E}=(b^{E}_{1},\ldots,b^{E}_{s})^{T},
\end{split}
\end{equation*}
and the s-stage IMEX RK method can be represented by a double Butcher tableau,
\begin{table}[!ht]
\begin{tabular}{c|c}
$c^{I}$ & $A^{I}$ \\
\hline & $(b^{I})^T$
\end{tabular},$\quad \quad$
\begin{tabular}{c|c}
$c^{E}$ & $A^{E}$ \\
\hline & $(b^{E})^T$
\end{tabular},\\
\end{table}\\
where $c^{I}_{i}=\sum\limits_{j=1}^{s} a^{I}_{ij},\,\,c^{E}_{i}=\sum\limits_{j=1}^{s} a^{E}_{ij},\,\,i=1,\ldots,s.$

We rewrite \eqref{RRK} into
\begin{equation}\label{DRRKr}
\begin{cases}
u^\varepsilon_{i}=u^{\varepsilon,n}_{\gamma}+\tau\sum\limits^{i}_{j=1}a^{I}_{ij}g^{I}_{j}+\tau\sum\limits^{i-1}_{j=1}a^{E}_{ij}g^{E}_{j},\,\, i=1,\ldots,s,\\
u^{\varepsilon,n+1}=u^{\varepsilon,n}_{\gamma}+\tau\sum\limits^{s}_{j=1}b^{I}_{j}g^{I}_{j}+\tau\sum\limits^{s}_{j=1}b^{E}_{j}g^{E}_{j},\\
u^{\varepsilon,n+1}_{\gamma}=u^{\varepsilon,n+1}+(\gamma_{n}-1)(u^{\varepsilon,n+1}-u^{\varepsilon,n}_{\gamma}).
\end{cases}
\end{equation}
\begin{thm}
For the semi-discrete scheme \eqref{RRK}, if we set
\begin{equation}\label{gammarrk}
\gamma_{n}=\begin{cases}
1,\quad &\sum\limits^{s}_{j=1}(b^{I}_{j}g^{I}_{j}+\tau\gamma_{n}\sum\limits^{s}_{j=1}b^{E}_{j}g^{E}_{j})=0, \\ \frac{-2\sum\limits^{s}_{j=1}\Re(b^{I}_{j}g^{I}_{j}+b^{E}_{j}g^{E}_{j},u^n_{\gamma})}{\tau\|\sum\limits^{s}_{j=1}b^{I}_{j}g^{I}_{j}+\sum\limits^{s}_{j=1}b^{E}_{j}g^{E}_{j}\|^2}, &\sum\limits^{s}_{j=1}(b^{I}_{j}g^{I}_{j}+\sum\limits^{s}_{j=1}b^{E}_{j}g^{E}_{j})\neq 0.
\end{cases}
\end{equation}
\eqref{RRK} would satisfy the conservation of mass:
\begin{equation}\label{massDRRK}
\|u^{\varepsilon,n+1}_{\gamma}\|^{2}=\|u^{\varepsilon,n}_{\gamma}\|^{2}.
\end{equation}
\end{thm}
\begin{proof}
From \eqref{RRK} and the definition of \eqref{gammarrk}, we have
\begin{equation}
\begin{split}
\|u^{\varepsilon,n+1}_{\gamma}\|^{2}-\|u^{\varepsilon,n}_{\gamma}\|^{2}&=\|u^{\varepsilon,n}_{\gamma}+\tau\gamma_{n}\sum\limits^{s}_{j=1}b^{I}_{j}g^{I}_{j}+\tau\gamma_{n}\sum\limits^{s}_{j=1}b^{E}_{j}g^{E}_{j}\|^2-\|u^{\varepsilon,n}_{\gamma}\|^{2} \\
&=\tau^2\gamma_{n}^2\|\sum\limits^{s}_{j=1}b^{I}_{j}g^{I}_{j}+\sum\limits^{s}_{j=1}b^{E}_{j}g^{E}_{j}\|^2+ 2\tau\gamma_{n}\sum\limits^{s}_{j=1}\Re(b^{I}_{j}g^{I}_{j}+b^{E}_{j}g^{E}_{j}, u^{\varepsilon,n}_{\gamma})\\
&=0.
\end{split}
\end{equation}
 This ends the proof.
\end{proof}
\begin{lemma}\label{gammaorder}\cite{Li2023Implicit}
Suppose that the given IMEX RK method is $p$-th order accurate with $p \geq 2$. For sufficiently small $\tau$, the relaxation coefficient $\gamma_n$ defined in \eqref{gammarrk} satisfies
\begin{equation}\label{gammaorder1}
\gamma_n=1+\mathcal{O}(\tau^{p-1}).
\end{equation}
\end{lemma}
\begin{thm}\label{truncationerr}
The truncation error of IMEX RRK method \eqref{DRRKr} is $O(\tau^{p+1})$.
\end{thm}
\begin{proof}
Since the convergence order of IMEX RK method is $O(\tau^p)$, putting the exact solution into the second equation of \eqref{DRRKr}, one can obtain
\begin{equation*}\label{truncationerreq}
u^\varepsilon(\hat{t}_n+\tau)=\Phi_n(u^\varepsilon(t))+O(\tau^{p+1}),
\end{equation*}
where $$\Phi_n(u^\varepsilon(t))=u^\varepsilon(\hat{t}_n)+\tau\sum\limits^{s}_{j=1}b^{I}_{j}g^{I}(u^\varepsilon(\hat{t}_{nj}))+\tau\sum\limits^{s}_{j=1}b^{E}_{j}g^{E}(u^\varepsilon(\hat{t}_{nj})),$$
and $\hat{t}_{nj}=\hat{t}_{n}+c_j \tau,\,j=1,\ldots,s$.

Then substitute the exact solution into the third equation of \eqref{DRRKr}, we obtain the truncation error of the third equation of \eqref{DRRKr} as follows
\begin{equation*}
\begin{split}
\mathcal{T}^{n+1}&=u^\varepsilon(\hat{t}_{n+1})-\Phi_n(u^\varepsilon(t))-(\gamma_{n}-1)\big(\Phi_n(u^\varepsilon(t))-u^\varepsilon(\hat{t}_{n})\big)\\
&=u^\varepsilon(\hat{t}_{n+1})-\big(u^\varepsilon(\hat{t}_n+\tau)+O(\tau^{p+1})\big)-(\gamma_{n}-1)\big(u^\varepsilon(\hat{t}_n+\tau)+O(\tau^{p+1})-u^\varepsilon(\hat{t}_{n})\big)\\
&=u^\varepsilon(\hat{t}_{n+1})-u^\varepsilon(\hat{t}_n+\tau)-(\gamma_{n}-1)\big(u^\varepsilon(\hat{t}_n+\tau)-u^\varepsilon(\hat{t}_{n})\big)+O(\tau^{p+1})+O\big((\gamma_{n}-1)\tau^{p+1}\big)\\
&=u^\varepsilon(\hat{t}_{n}+\tau+(\gamma_{n}-1)\tau)-u^\varepsilon(\hat{t}_n+\tau)-(\gamma_{n}-1){u^{\varepsilon}}^{\prime}(\hat{t}_n+\tau)\tau+O(\tau^{p+1})\\
&=u^\varepsilon(\hat{t}_{n}+\tau)+(\gamma_{n}-1)\tau {u^{\varepsilon}}^{\prime}(\hat{t}_n+\tau)+O\big((\gamma_{n}-1)^2\tau^{2}\big)\\
&\quad -u(\hat{t}_n+\tau)-(\gamma_{n}-1)\tau {u^{\varepsilon}}^{\prime}(\hat{t}_n+\tau)+O(\tau^{p+1})=O(\tau^{p+1}),
\end{split}
\end{equation*}
where we used $\hat{t}_{n+1}=\hat{t}_{n}+\gamma_{n}\tau$.
This ends the proof.
\end{proof}

\section{Numerical Experiments}

In the numerical test, the following RK methods will be used in time:
\begin{enumerate}
  \item[1.] RK(1,2), Implicit-explicit midpoint \cite{Herty2013Implicit};
  \item[2.] RK(2,3), two-stage, third-order RK \cite{Herty2013Implicit};
  \item[3.] RK(6,4), also named ARK4(3)6L[2]SA in \cite{Kennedy2003Additive};
  \item[4.] RK(8,5), also named ARK5(4)8L[2]SA in \cite{Kennedy2003Additive},
\end{enumerate}
and we apply Fourier pseudo-spectral method in space where can be used FFT method. We only show
numerical results for the ERLogSE, the results of ERSSE are quite similar which are omitted for brevity.
\subsection{Accuracy test}
Consider the LogSE \eqref{logse} with the exact Gaussian solution \cite{Carles2018Universal} in $1$-dimension as follows,
\begin{equation}\label{2dGex}
u(x,t)=b\exp\big(\ri\big(x\cdot \varsigma -(a+| \varsigma|^2)t\big)+(\lambda/2)| x-2 \varsigma t|^2 \big), \quad  x \in \mathbb{R}, \quad t \geq 0,
\end{equation}
where $a=-\lambda (1-\ln|b|^2)$, $b\in \mathbb{R}$ and $\lambda, \,\varsigma\in \mathbb{R}$ are given constants.
To quantify the numerical errors, we define the following error functions:
\begin{equation}
\begin{split}
e^{\varepsilon}(t_n)&=u^{\varepsilon}_{ref}(\cdot,t_n)- u^{\varepsilon,n},\quad \hat{e}^{\varepsilon}(t_n)=u_{ex}(\cdot,t_n)-u^{\varepsilon}_{ref}(\cdot,t_n),\\
\hat{e}^{\varepsilon}_{\rho}(t_n)&=\rho_{ex}(\cdot,t_n)- \rho^{\varepsilon}_{ref}(\cdot,t_n)=|u_{ex}(\cdot,t_n)|^2-|u^{\varepsilon}_{ref}(\cdot,t_n)|^2,\quad e^{\varepsilon}_{E}=E(u_0)- E^{\varepsilon}_k (u_0).
\end{split}
\end{equation}
Here $u^{\varepsilon,n}_j$ is the numerical solution, $u_{ex}$ is the exact solution of LogSE \eqref{logse}, $u^{\varepsilon}_{ref}$ is the `exact' solution of ERLogSE \eqref{rlogse} with $h_e=\frac{5}{2^{11}},\,\tau_e=10^{-5}$. The initial data is taken as $u_0=u(x,0)$ in \eqref{2dGex} and the  boundary conditions are given such that the exact solution. Here $\varsigma=0,\,T=1,\, \Omega=[-10,10]$, $\lambda=-1, b=1$.
\subsubsection{Convergence rate of the regularized model}
\begin{figure}[!ht]
  \centering
  \subfigure{\includegraphics[width=0.3\textwidth]{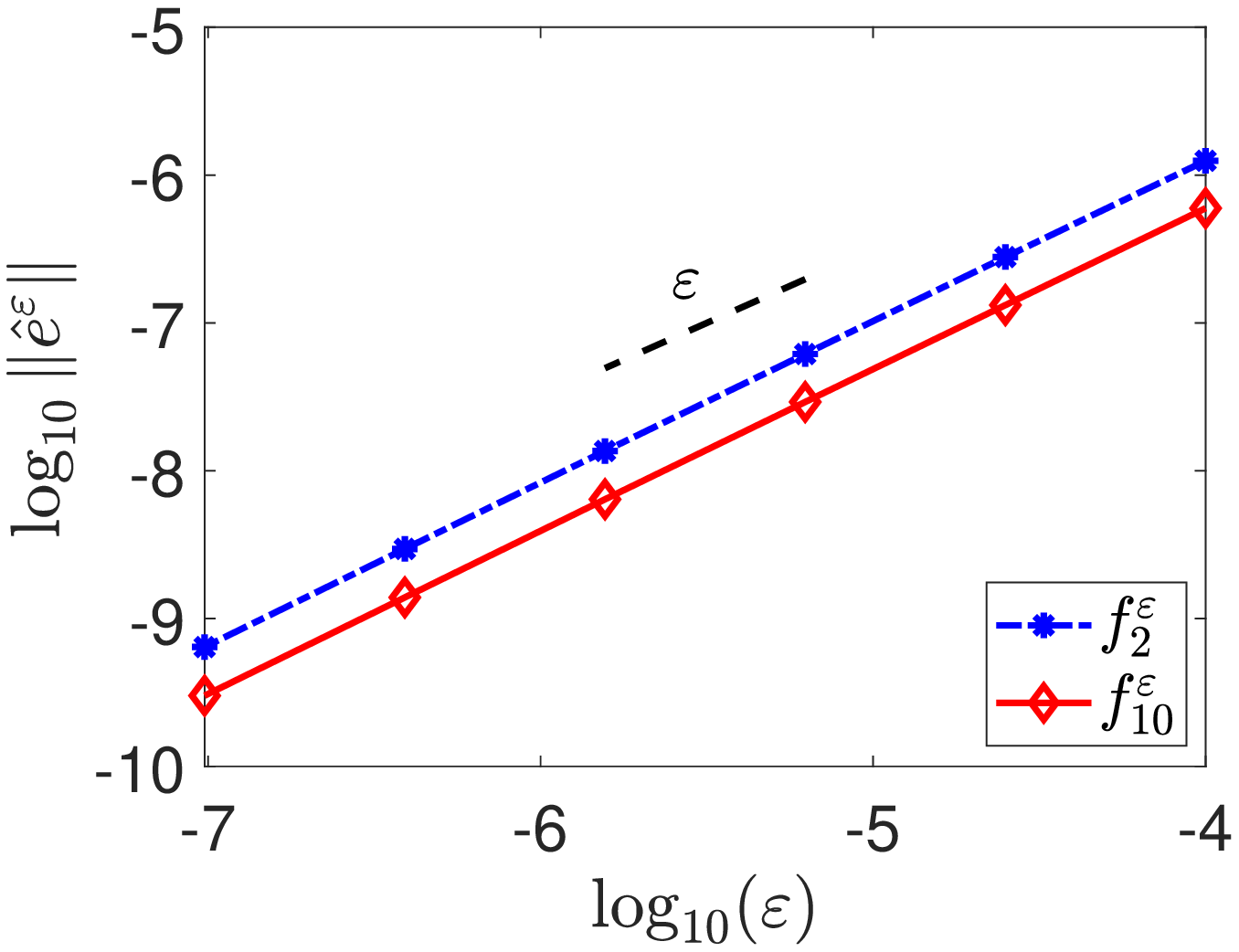}}
  \subfigure{\includegraphics[width=0.3\textwidth]{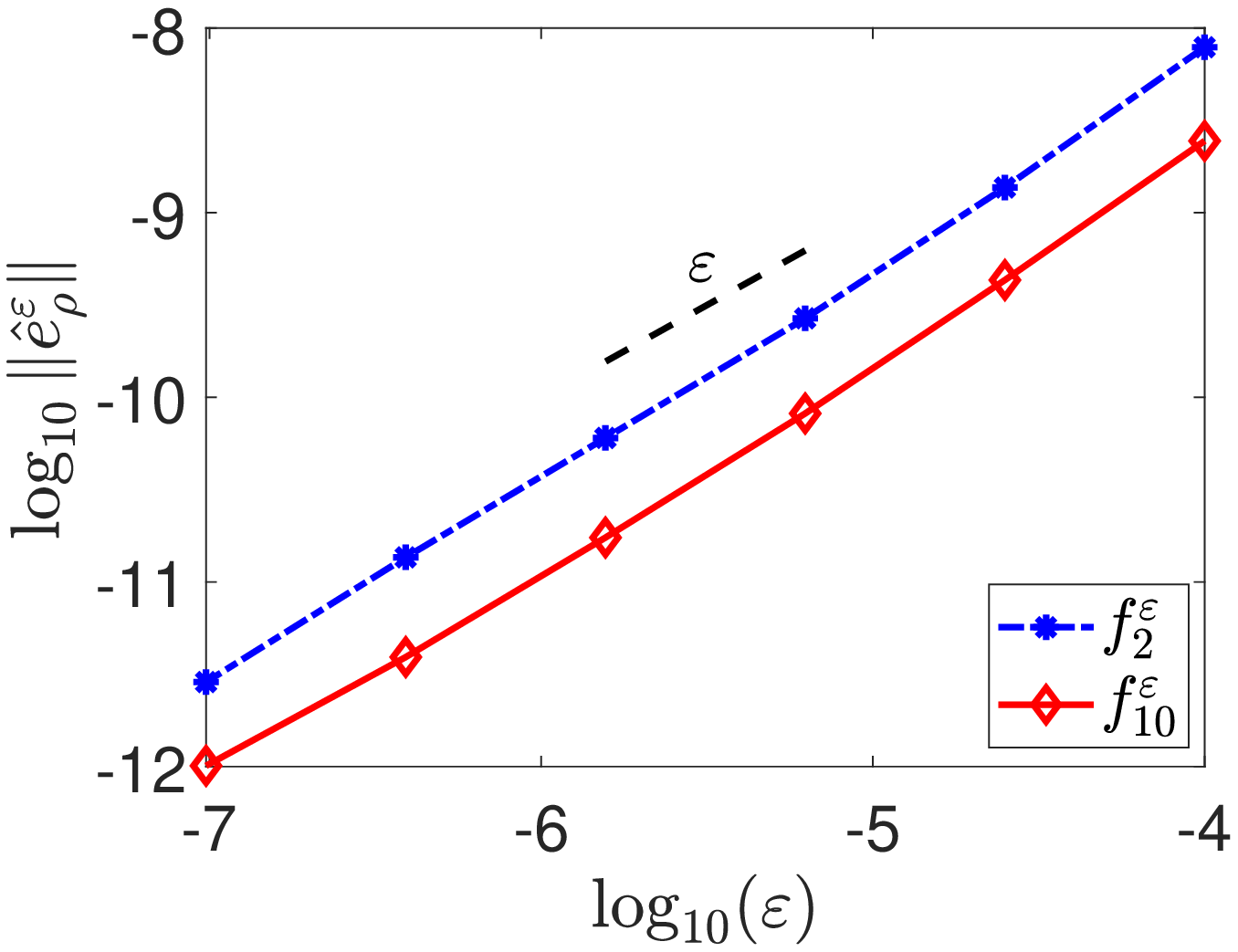}}
  \subfigure{\includegraphics[width=0.3\textwidth]{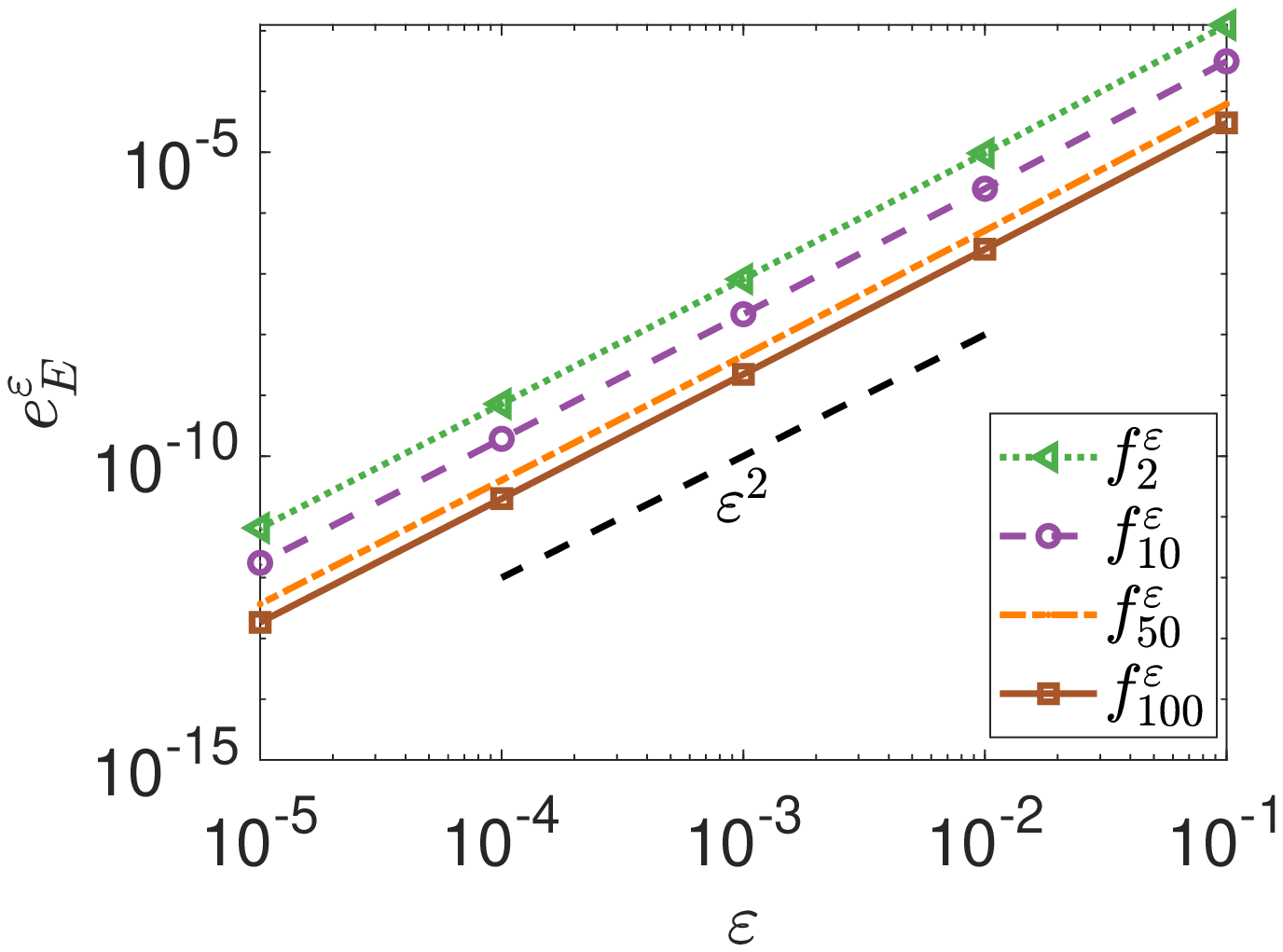}}
\caption{Convergence order ERLogSE to LogSE, i.e. $\|\hat{e}^{\varepsilon}(t=1)\|,\,\|\hat{e}^{\varepsilon}_{\rho}(t=1)\|,\,e^{\varepsilon}_{E}(t=1)$.}\label{ehatorder}
\end{figure}

Here we consider the error between the solutions of the ERLogSE \eqref{rlogse} and the LogSE \eqref{logse}. Fig. \ref{ehatorder} depicts $\|\hat{e}^{\varepsilon}\|,\,\|\hat{e}^{\varepsilon}_{\rho}\|,\,e^{\varepsilon}_{E}$ at $t=1$ with different regularized nonlinearities $f^{\varepsilon}_2,\,f^{\varepsilon}_{10}$ computed by IMEX RRK(2,3), and other RRK schemes are not depicted here for the sake of conciseness. We can draw the following conclusions: (i) The solution of ERLogSE \eqref{rlogse} converges linearly to LogSE \eqref{logse} in terms of $\varepsilon$; (ii) The density $\rho^{\varepsilon}$ of the solution of ERLogSE \eqref{rlogse} exhibits linear convergence to  that of the LogSE \eqref{logse} in terms of $\varepsilon$; (iii) The regularized energy $E^{\varepsilon}_k$ converges quadratically to the original energy $E$; (iv)  Across these three figures, it is evident that for any given $\varepsilon$, the solutions with larger values of $n$ yield superior performance.
\subsubsection{Convergence rate of the numerical scheme}
\begin{figure}[!ht]
  \centering
\subfigure[$n=2,\,\varepsilon=10^{-5}$]{\includegraphics[width=0.23\textwidth]{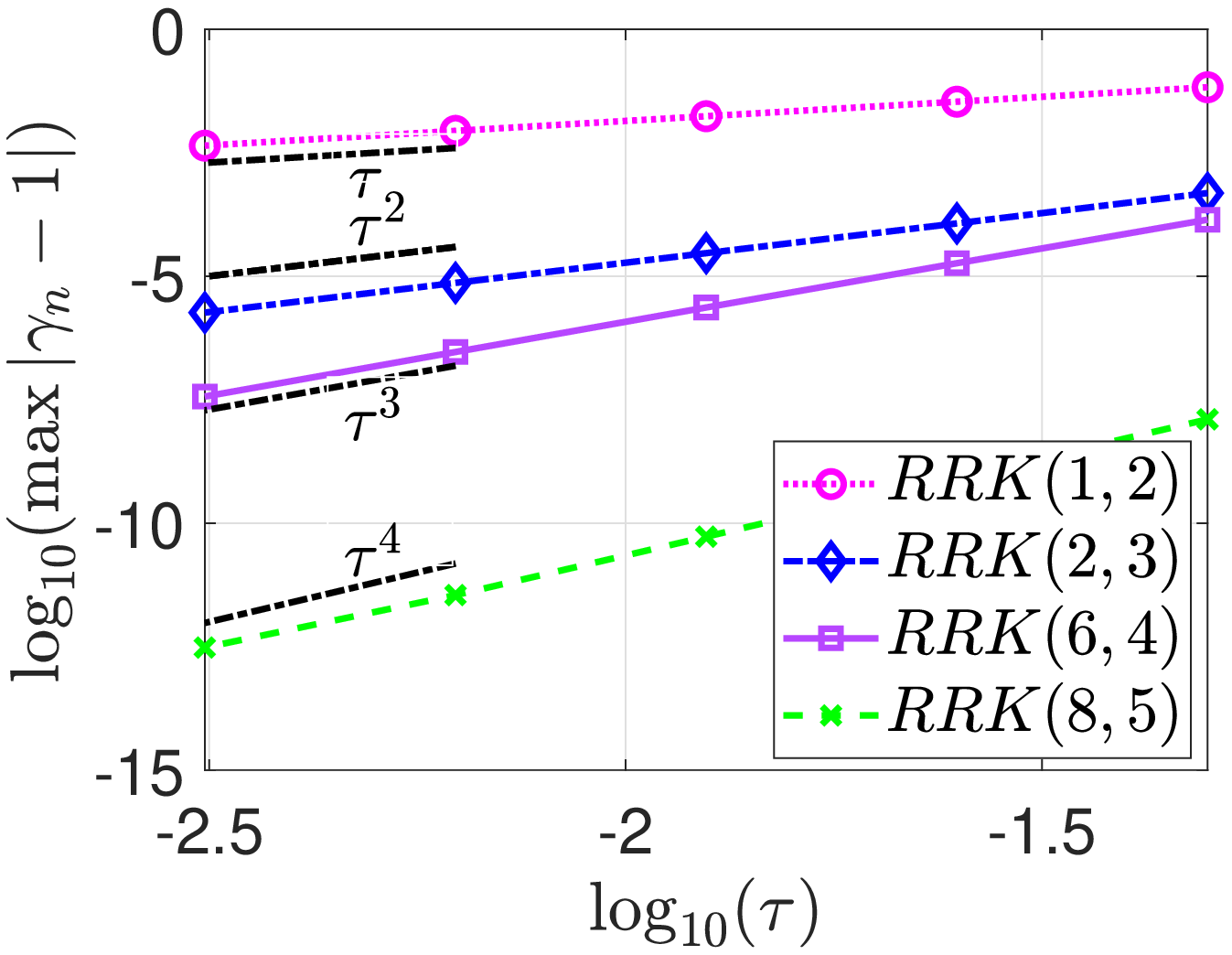}}
  \subfigure[$n=2,\,\varepsilon=10^{-10}$]{\includegraphics[width=0.23\textwidth]{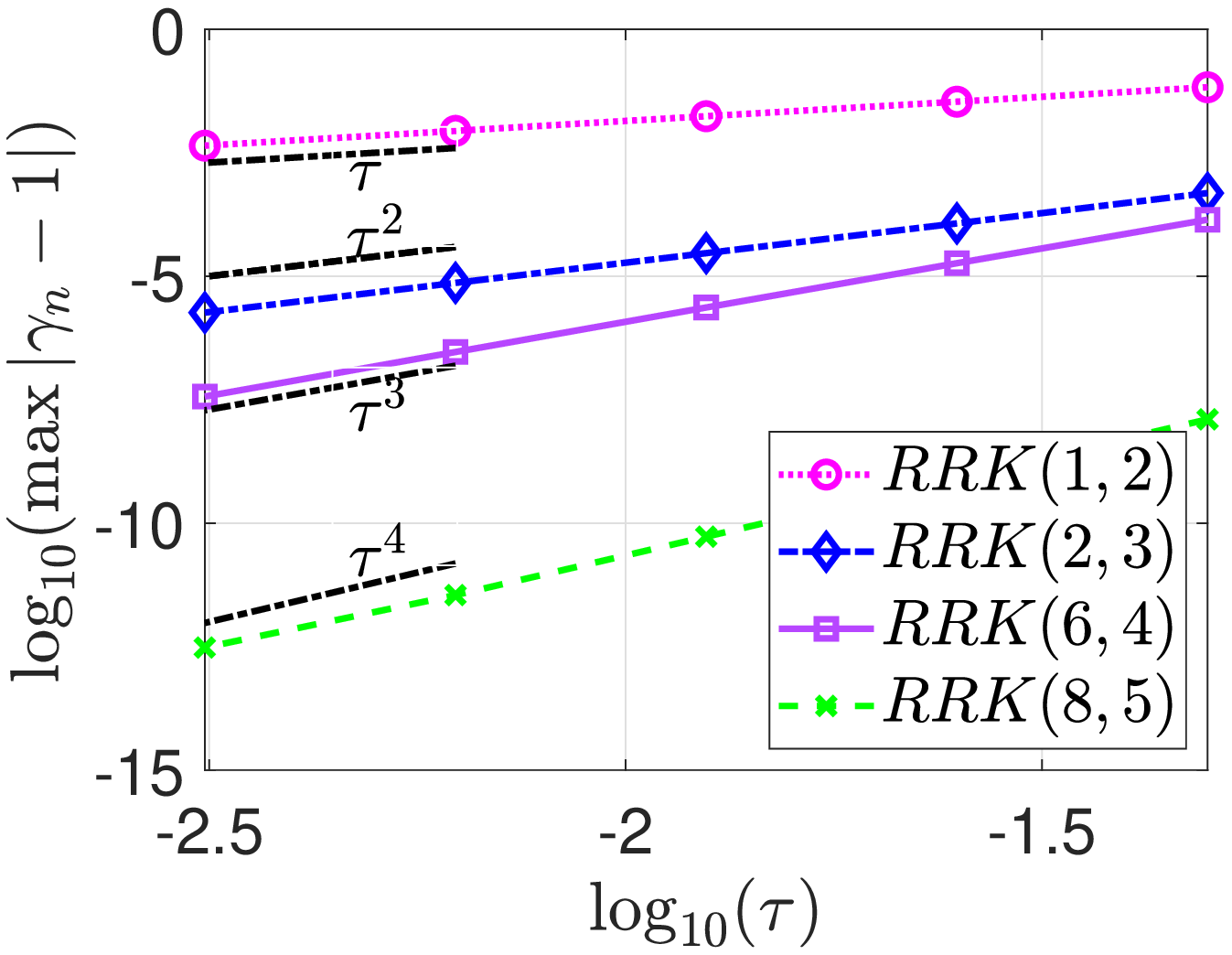}}
  \subfigure[$n=4,\,\varepsilon=10^{-5}$]{\includegraphics[width=0.23\textwidth]{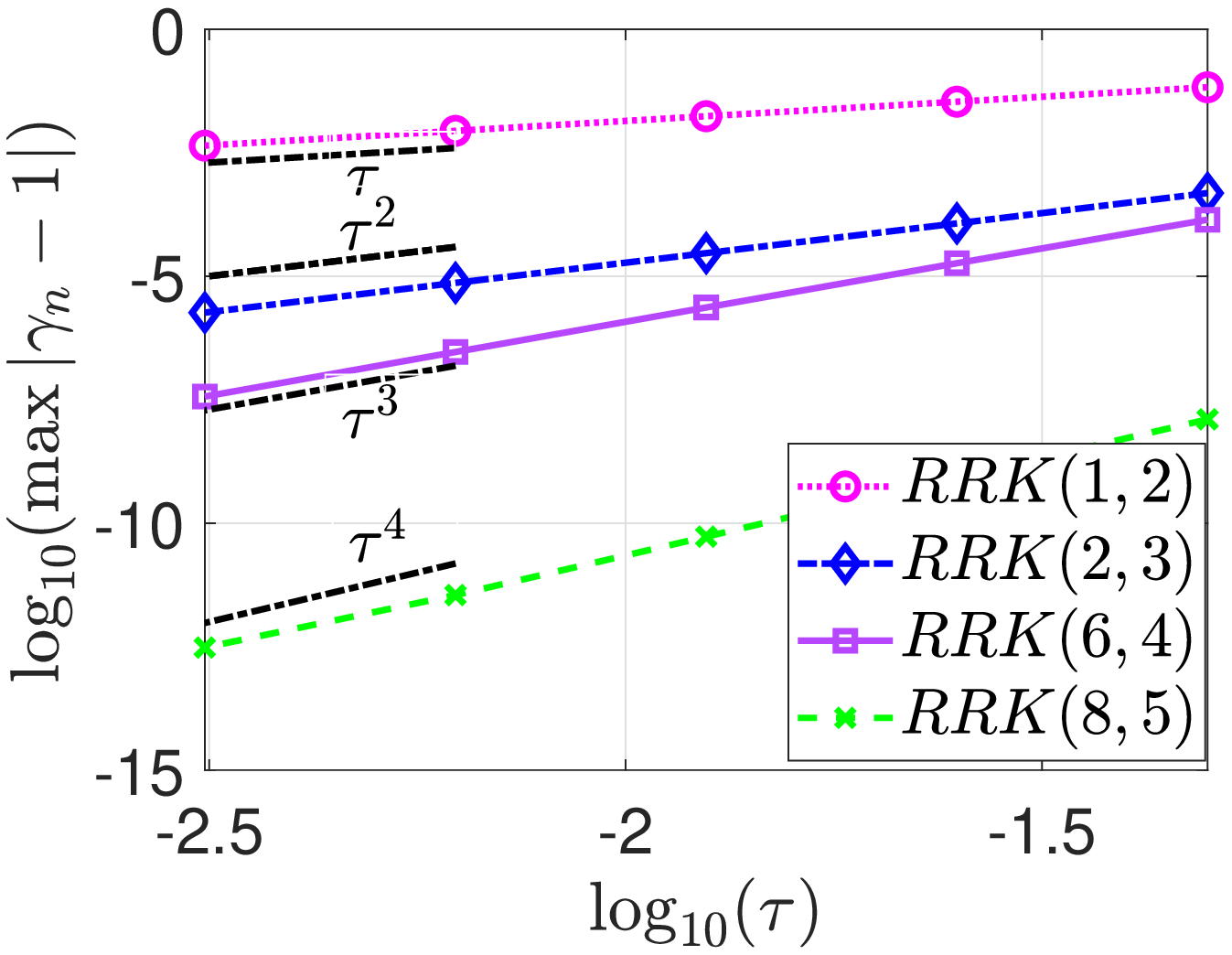}}
  \subfigure[$n=4,\,\varepsilon=10^{-10}$]{\includegraphics[width=0.23\textwidth]{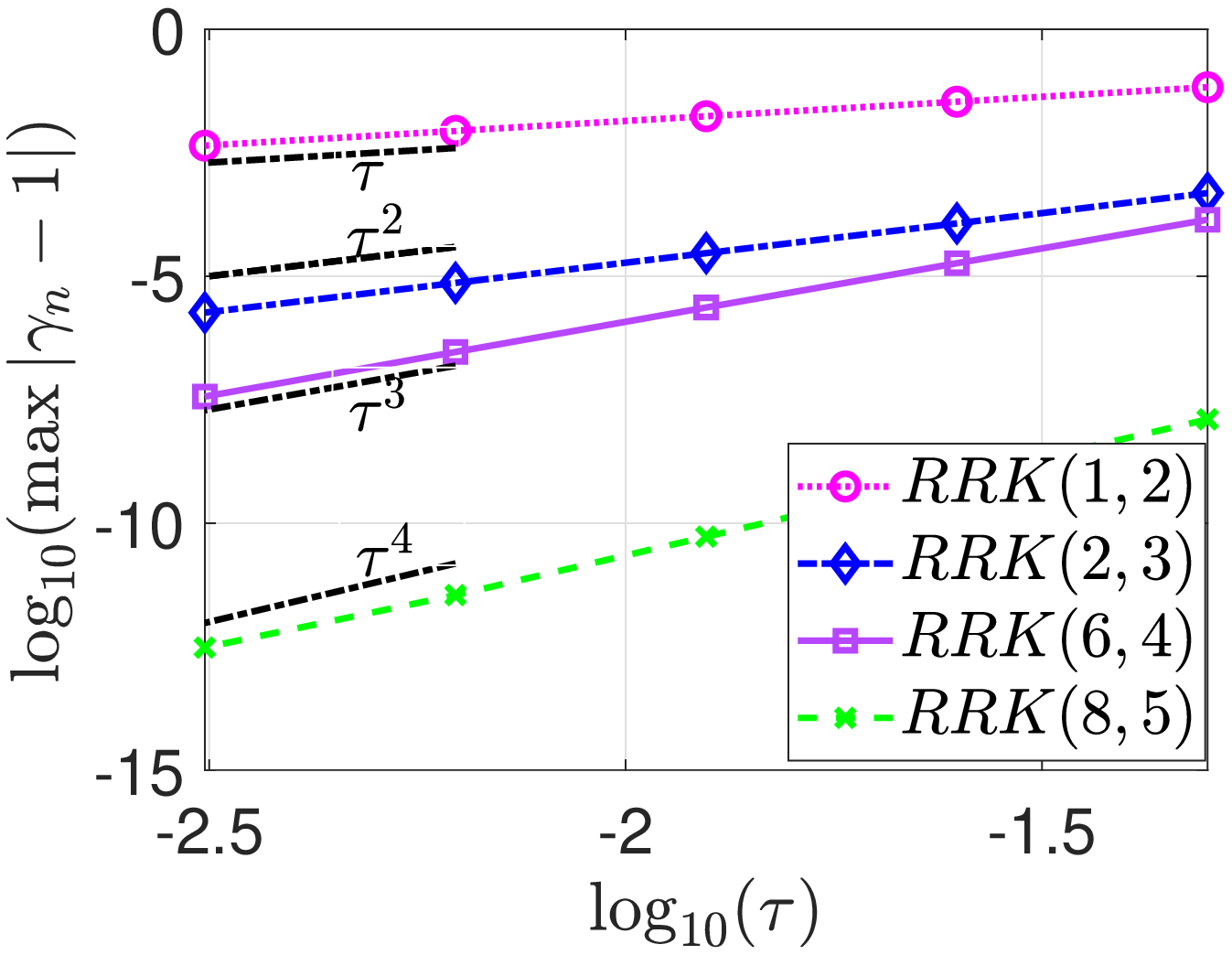}}
\caption{Convergence order of $max(|\gamma_n-1|)$ for different IMEX RRK schemes.}\label{gammaorderfig}
\end{figure}

Firstly, we investigate the relaxation coefficient $\gamma_n$ at $t=1$ for various values of $n,\,\varepsilon$ and the results shown in Fig.\ref{gammaorderfig}, substantiate Lemma \ref{gammaorder}.

\begin{figure}[!ht]
  \centering
  \subfigure[RRK(1,2)]{\includegraphics[width=0.23\textwidth]{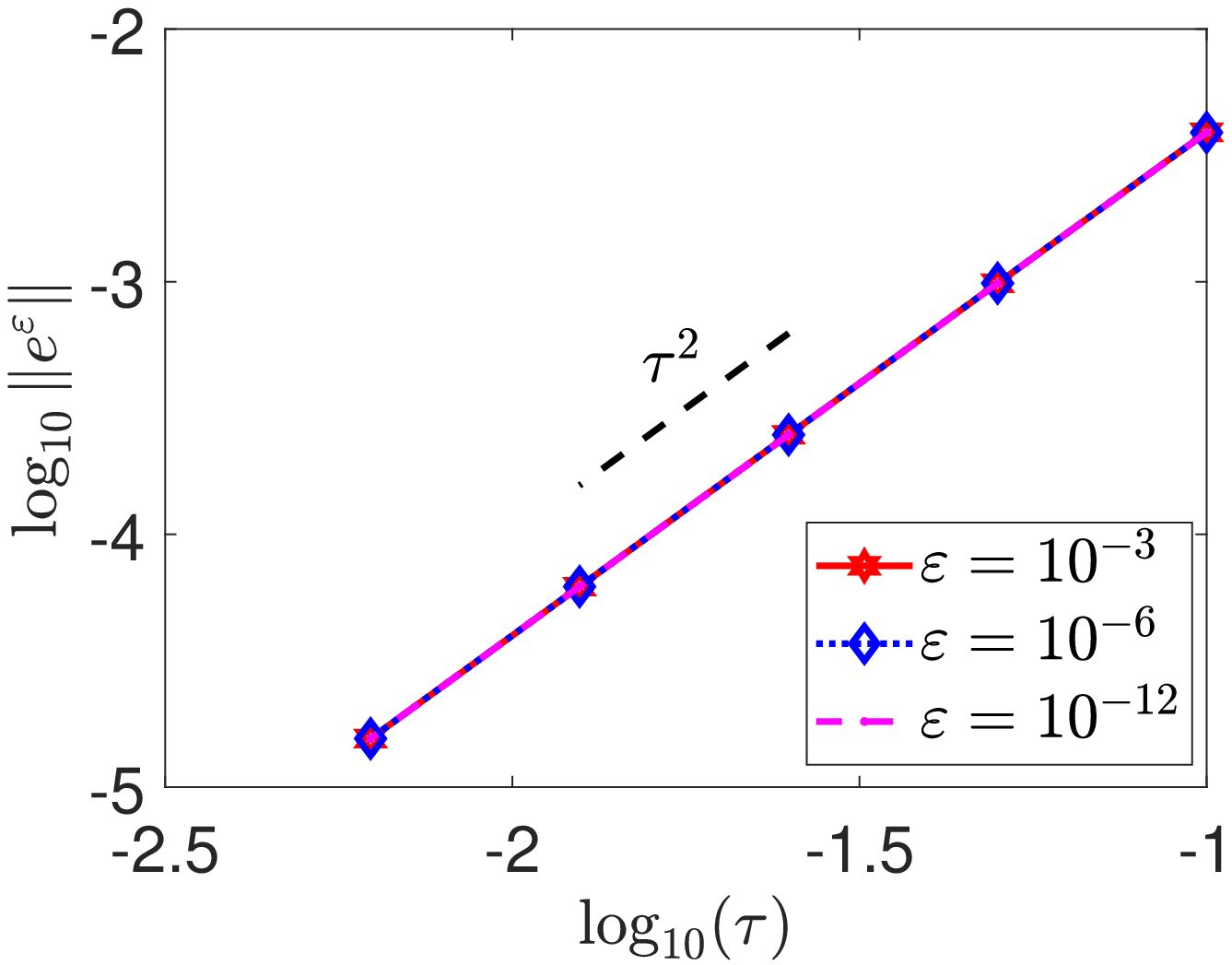}}
  \subfigure[RRK(2,3)]{\includegraphics[width=0.23\textwidth]{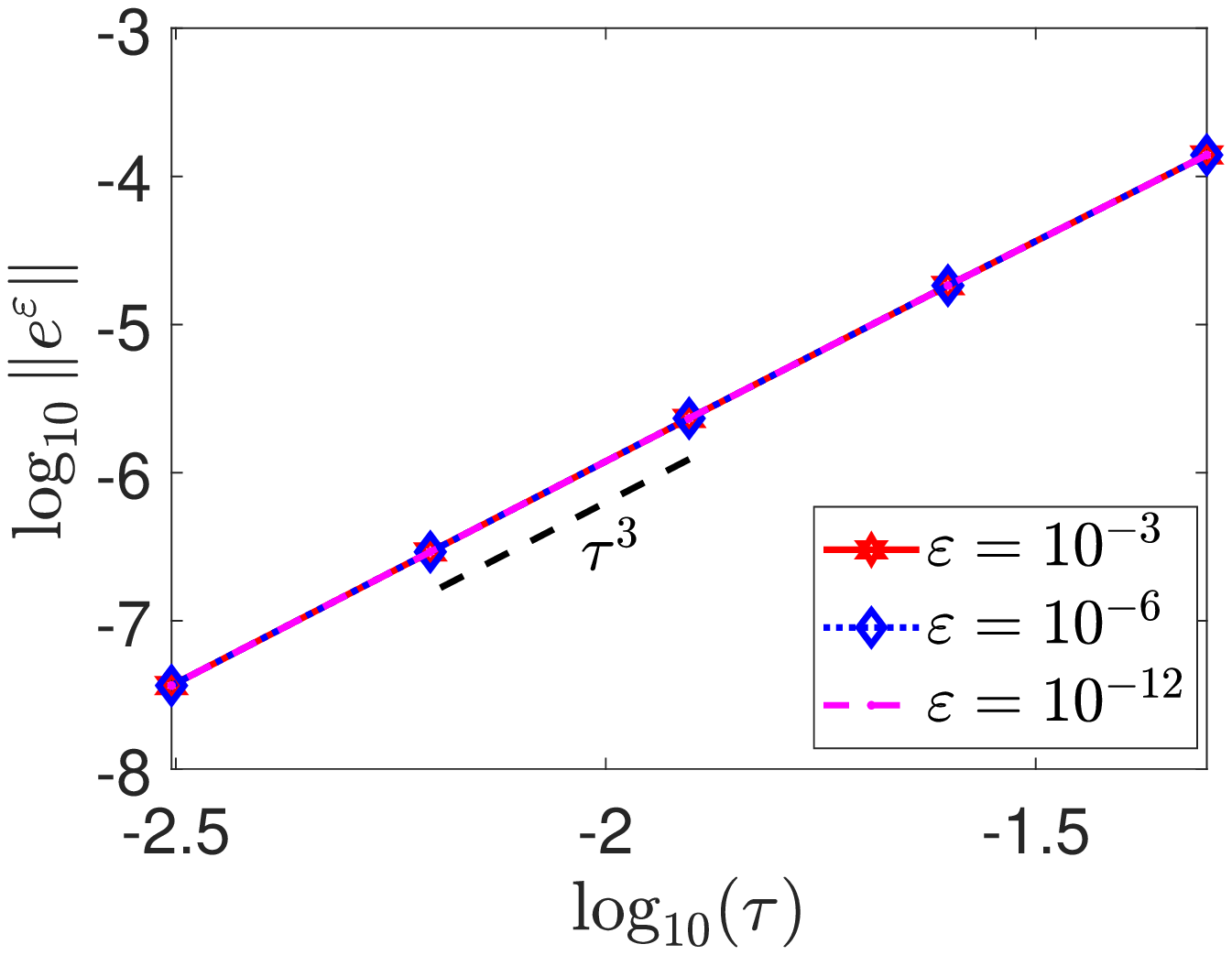}}
  \subfigure[RRK(6,4)]{\includegraphics[width=0.23\textwidth]{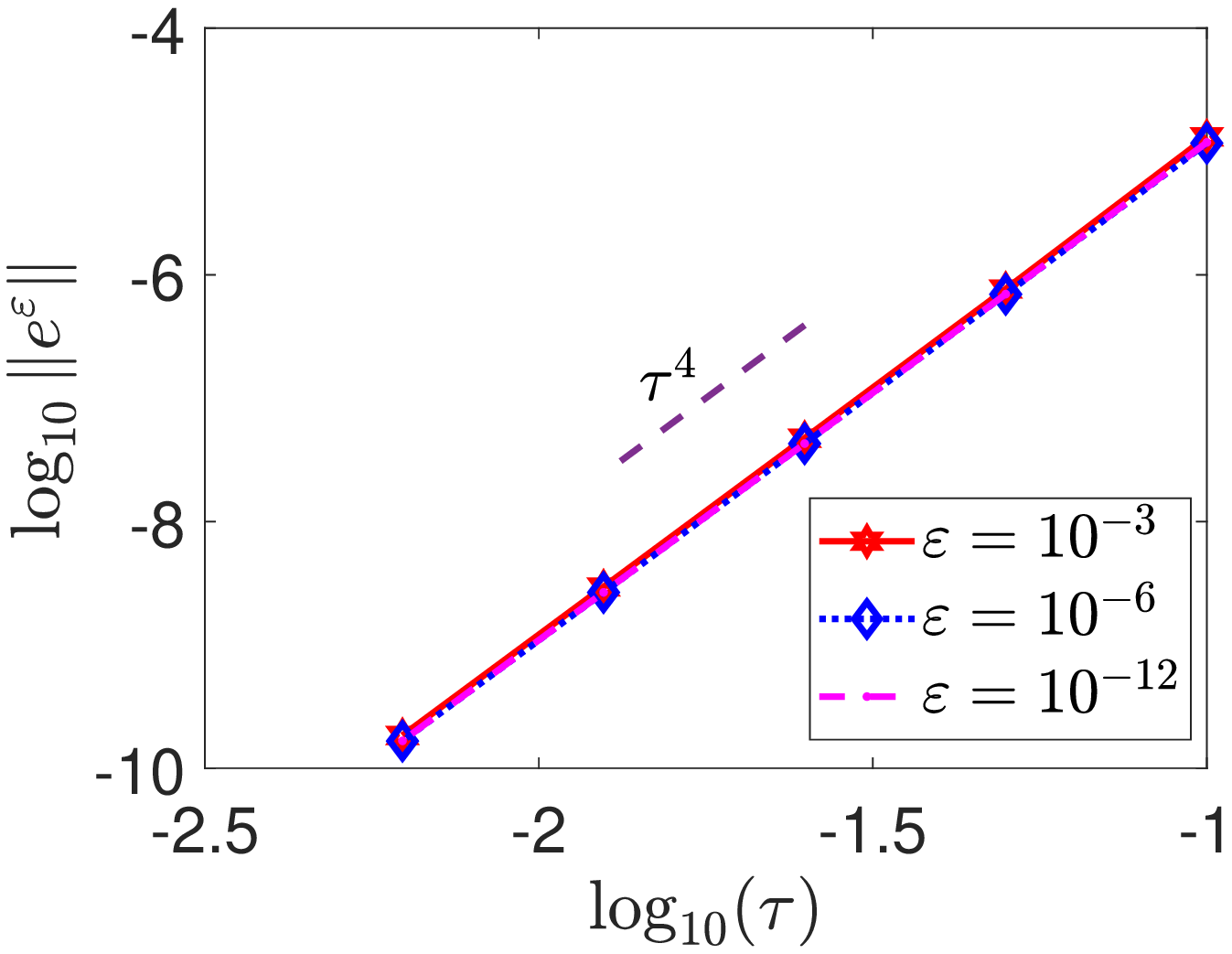}}
  \subfigure[RRK(8,5)]{\includegraphics[width=0.23\textwidth]{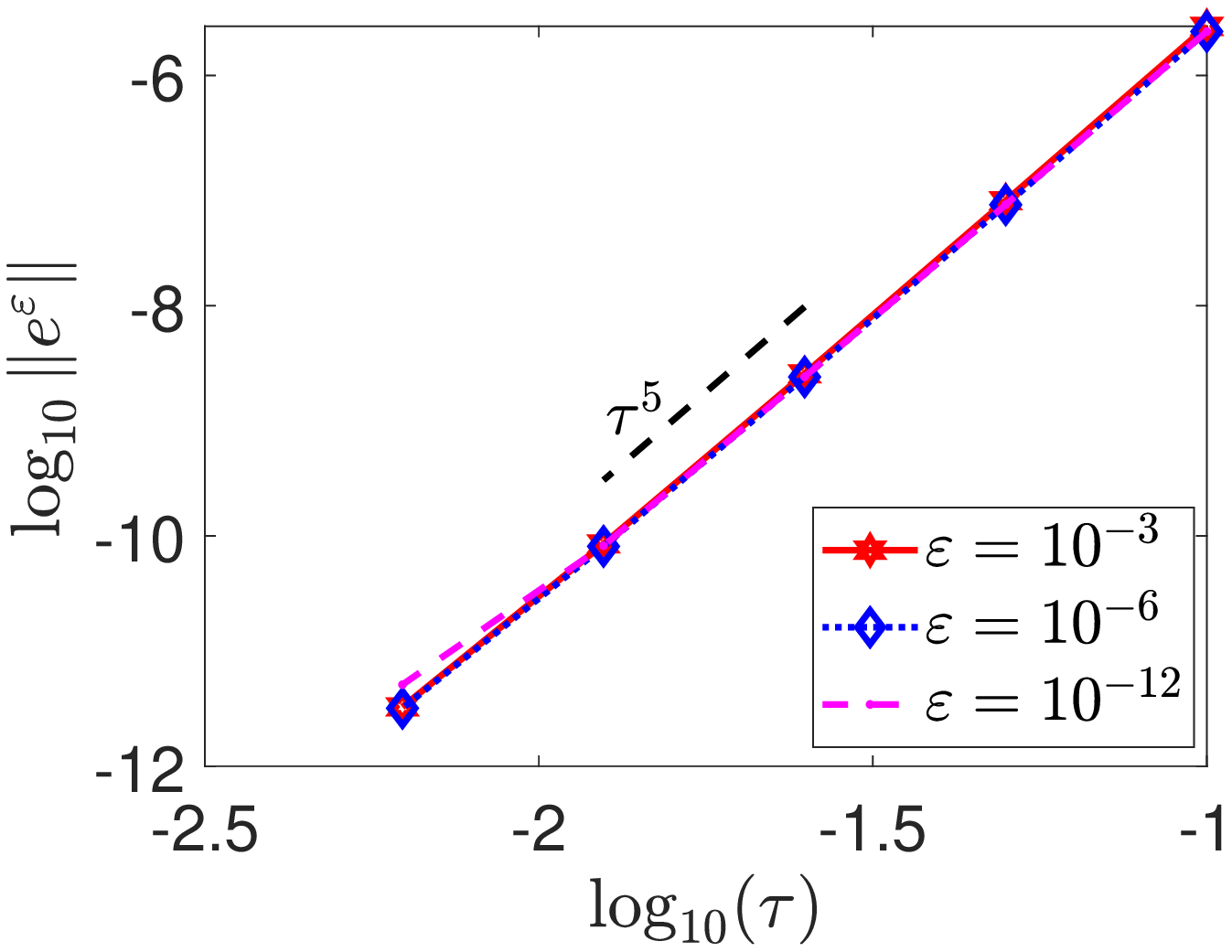}}
\caption{Convergence order in time of RRK with different $\varepsilon$.}\label{FIGRRK}
\end{figure}

Next, to assess the temporal convergence rate, we set the mesh size to $h=h_e$, and varied the time step as $\tau=2^{-j}\times 10^{-1}$ for $j=1,2,\ldots,5$ with $n=2,\,n=4$ under different $\varepsilon$. Fig.\ref{FIGRRK} displays the temporal convergence rates of the IMEX RRK methods, calculated using various schemes. The RRK methods show a convergence rate of $p$, confirming the theoretical predictions in Theorem \ref{truncationerr}.

Additionally, to evaluate the spatial convergence rate of the Fourier pseudo-spectral method, we fixed the time step
$\tau=\tau_e$, and $\varepsilon=10^{-6}$ and selected different numbers of grid points $N_{j}=8+2(j+1),j=1,2,\ldots,5$. Table \ref{spaceorder} illustrates the spatial convergence rates of the IMEX RRK method, with the RRK(2,3) scheme being particularly noteworthy.

\begin{table}
\centering
\caption{Space accuracy tests of the Fourier spectral discretization ($\tau=1\times 10^{-5},\,\varepsilon=10^{-6},\,T=1$).}
\begin{tabular}{c|c|cc}
\hline \multirow{2}{*}{\text { RK (Stage,\,Order) }} & & \multicolumn{2}{c}{\text { RRK }}  \\
& $N$ & $e^{\varepsilon}$ &  \text { order }  \\
\hline \multirow{3}{*}{\text { (2,3) }}
& 12 & $2.54 \mathrm{e}-01$ & -    \\
& 14 & $1.38 \mathrm{e}-01$ & 3.95   \\
& 16 & $5.59 \mathrm{e}-02$ & 6.77  \\
& 18 & $1.66 \mathrm{e}-02$ & 10.30  \\
& 20 & $3.73 \mathrm{e}-03$ & 14.17  \\
\hline
\end{tabular}\label{spaceorder}
\end{table}

\subsection{Dynamics}
In this section, we investigate long time dynamics of LogSE with Gaussian-type initial datum in 1D. To this end, we fix $\lambda =-1, x\in [-40,40],\tau=2\times10^{-3},h=\frac{5}{16}$. The initial datum is chosen as
\begin{equation}\label{Gaini}
u_{0}(x)=\sum_{k=1}^{2} b_{k} e^{-\frac{a_{k}}{2}\left(x-x_{k}\right)^{2}+\ri v_{k} x}, \quad x \in \mathbb{R},
\end{equation}
where $b_k,a_k, x_k$ and $v_k$ are real constants, i.e, the initial data is the sum of $2$ Gaussons  with velocity $v_k$ and initial location $x_k$.
We take $v_{1}=-v_{2}=2, x_{1}=-x_{2}=-30, b_{k}=a_{k}=1(k=1,2)$.

\begin{figure}[!ht]
  \centering
   \subfigure{\includegraphics[width=0.3\textwidth]{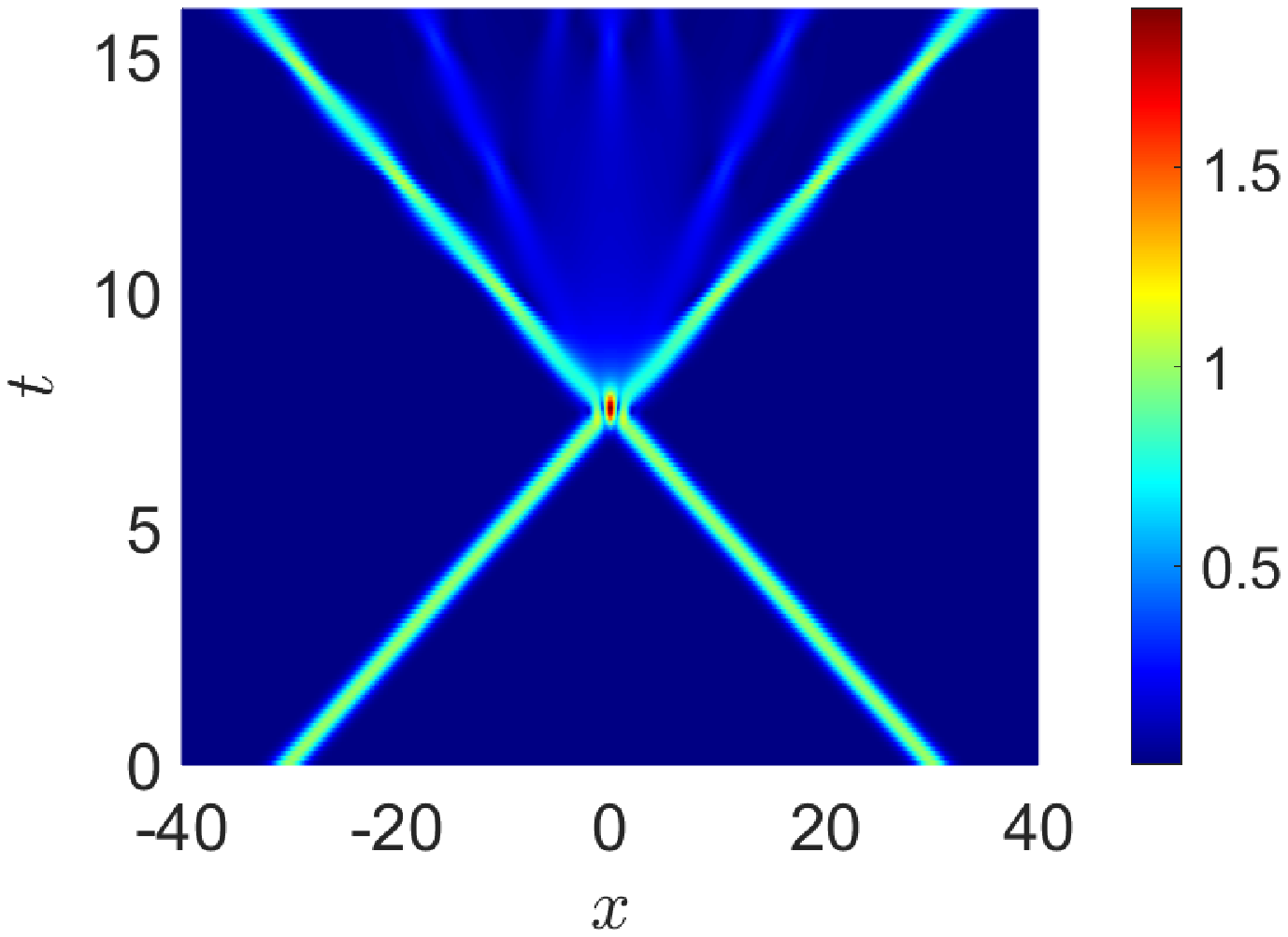}}
   \subfigure{\includegraphics[width=0.3\textwidth]{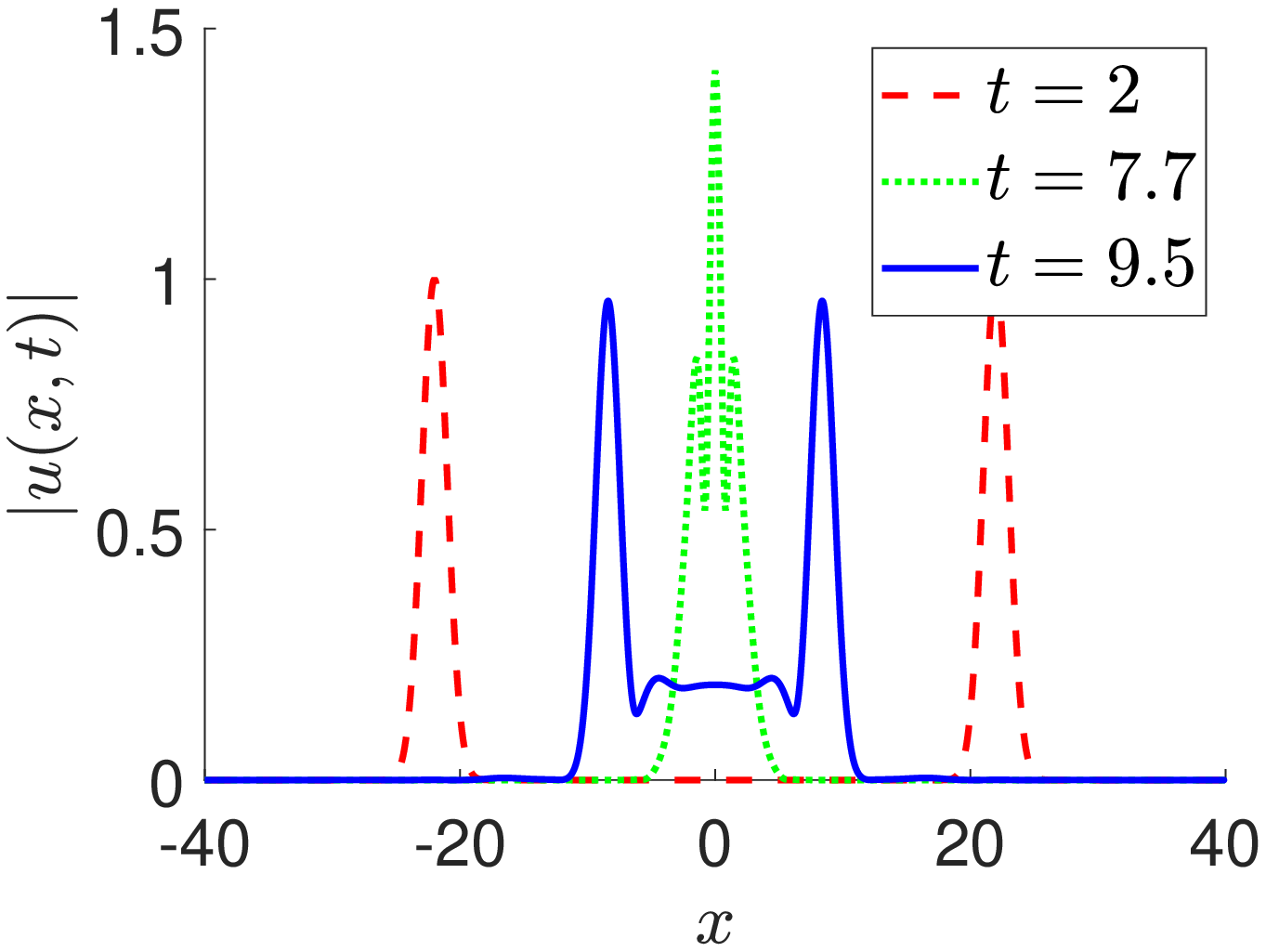}}
   \subfigure{\includegraphics[width=0.3\textwidth]{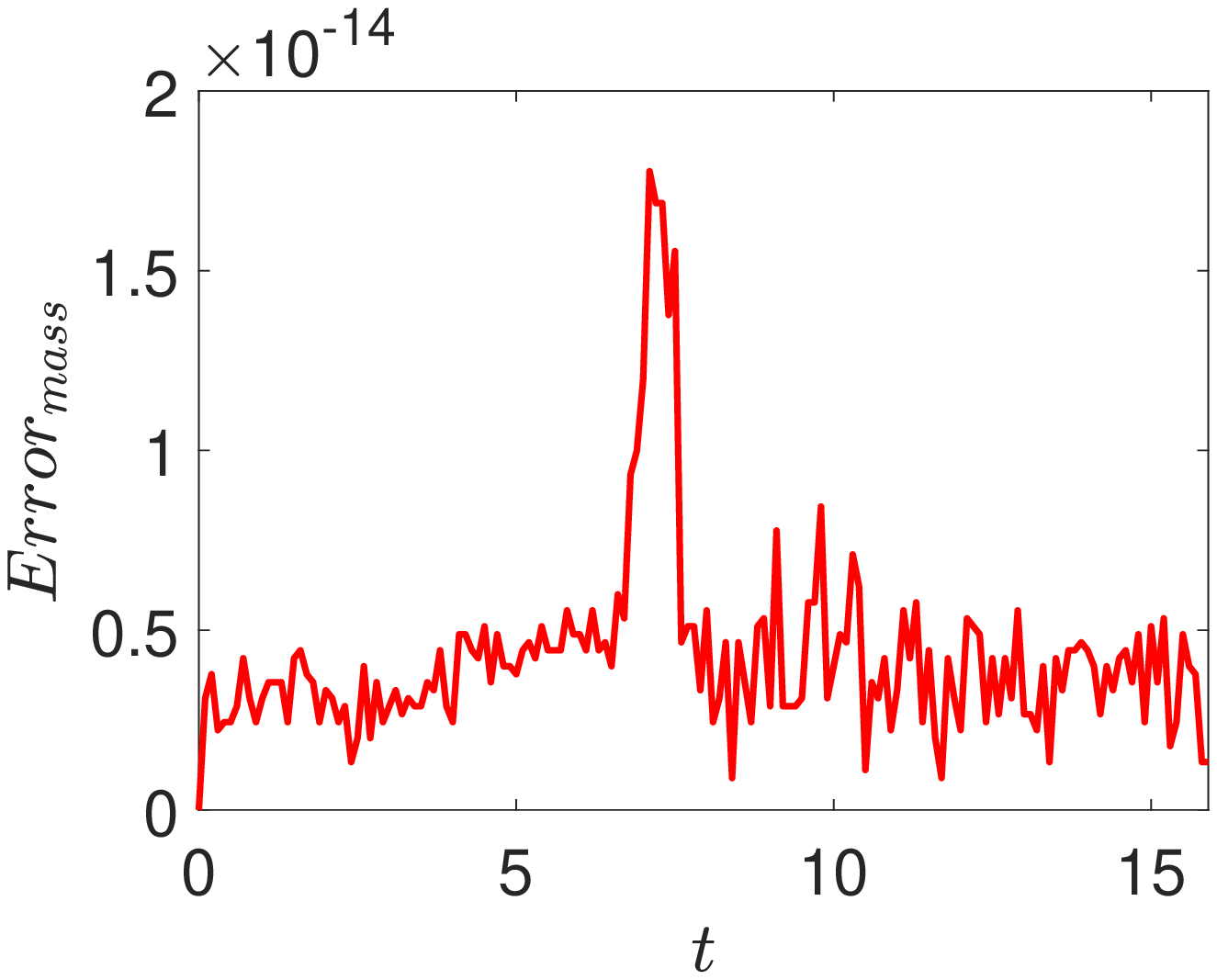}}
    \caption{Plots of $|u^{\varepsilon}(x,t)|$ (first column); $|u^{\varepsilon}(x,t)|$ at different time (second column) and evolution of mass error (third column) (Case for RRK(2,3)).}\label{Dynamtics}
\end{figure}
Analyzing Fig. \ref{Dynamtics}, we observe the following behaviors: (i) The mass is well conserved throughout the simulations, indicating the effectiveness of our numerical methods in preserving this fundamental property; (ii) In the case of moving Gaussons, they pass completely through each other, eventually moving separately. Oscillations occur during the interaction. 
\section{Conclusion}
We introduced an energy regularization approach to address the singularity present in the LogSE \eqref{logse} by employing a polynomial approximation. Subsequently, we developed a family of IMEX RRK methods, combined with the Fourier spectral method, to solve the ERLogSE \eqref{rlogse}. Our analysis demonstrated that these IMEX RRK methods not only inherit the mass conservation properties but also maintain the same convergence order as their standard counterparts. Furthermore, we conducted numerical experiments that substantiated our theoretical findings
\section*{Acknowledgements}
The work is supported by the National Natural Science Foundation of China (No. 12271523), Natural Science Foundation of Zhejiang Province (No. Q24A010031), 
Scientific Research Fund of Zhejiang Provincial Education Department (No. Y202352579), Natural Science Foundation of Jiangsu Province, (No. BK20240832), Defense Science Foundation of China (No. 2021-JCJQ-JJ-0538), Science \& Technology Innovation Program of Hunan Province (No. 2022RC1192).

\section*{References}
\bibliographystyle{plain}
\bibliography{Full}

\end{document}